\documentclass[10pt]{amsart}
\usepackage{amsmath}
\usepackage{amsxtra}
\usepackage{amscd}
\usepackage{amsthm}
\usepackage{amsfonts}
\usepackage{amssymb}
\usepackage{eucal}
\usepackage[all]{xy}
\usepackage{graphicx}

\newtheorem{prop}[subsubsection]{Proposition}

\newtheorem{conj}[subsubsection]{Conjecture}
\newtheorem{thm}[subsubsection]{Theorem}

\theoremstyle{definition}

\theoremstyle{remark}
\newtheorem{rem}[subsubsection]{Remark}

\newcommand{\thmref}[1]{Theorem~\ref{#1}}

\newcommand{\secref}[1]{Sect.~\ref{#1}}

\newcommand{\propref}[1]{Proposition~\ref{#1}}

\newcommand{\conjref}[1]{Conjecture~\ref{#1}}

\numberwithin{equation}{section}

\newcommand{\nc}{\newcommand}
\nc{\renc}{\renewcommand}
\nc{\ssec}{\subsection}
\nc{\sssec}{\subsubsection}
\nc{\on}{\operatorname}

\nc\ol{\overline}
\nc\wt{\widetilde}
\nc\tboxtimes{\wt{\boxtimes}}
\nc\tstar{\wt{\star}}
\nc{\alp}{\alpha}

\nc{\ZZ}{{\mathbb Z}}
\nc{\NN}{{\mathbb N}}
\nc{\OO}{{\mathbb O}}
\renc{\SS}{{\mathbb S}}
\nc{\DD}{{\mathbb D}}
\nc{\GG}{{\mathbb G}}

\nc{\Fq}{{\mathbb F}_q}
\nc{\Fqb}{\ol{{\mathbb F}_q}}
\nc{\Ql}{\ol{{\mathbb Q}_\ell}}
\nc{\id}{\text{id}}
\nc\X{\mathcal X}

\nc{\Hom}{\on{Hom}}
\nc{\Lie}{\on{Lie}}
\nc{\Loc}{\on{Loc}_{\cG,\on{spec}}}
\nc{\Pic}{\on{Pic}(X)}
\nc{\Bun}{\on{Bun}}
\nc{\IC}{\on{IC}}
\nc{\Aut}{\on{Aut}}
\nc{\rk}{\on{rk}}
\nc{\Sh}{\on{Sh}}
\nc{\Perv}{\on{Perv}}
\nc{\pos}{{\on{pos}}}
\nc{\Conv}{\on{Conv}}
\nc{\Sph}{\on{Sph}}
\nc{\Sym}{\on{Sym}}
\nc{\BunBb}{\overline{\Bun}_B}
\nc{\BunNb}{\overline{\Bun}_N}
\nc{\BunTb}{\overline{\Bun}_T}
\nc{\BunBbm}{\overline{\Bun}_{B^-}}
\nc{\BunBbel}{\overline{\Bun}_{B,el}}
\nc{\BunBbmel}{\overline{\Bun}_{B^-,el}}
\nc{\Buno}{\overset{o}{\Bun}}
\nc{\BunPb}{{\overline{\Bun}_P}}
\nc{\BunBM}{\Bun_{B(M)}}
\nc{\BunBMb}{\overline{\Bun}_{B(M)}}
\nc{\BunPbw}{{\widetilde{\Bun}_P}}
\nc{\BunBP}{\widetilde{\Bun}_{B,P}}
\nc{\GUb}{\overline{G/U}}
\nc{\GUPb}{\overline{G/U(P)}}

\nc{\Hhom}{\underline{\on{Hom}}}
\nc\syminfty{\on{Sym}^{\infty}}
\nc\lal{\ol{\lambda}}
\nc\xl{\ol{x}}
\nc\thl{\ol{\theta}}
\nc\nul{\ol{\nu}}
\nc\mul{\ol{\mu}}
\nc\Sum\Sigma
\nc{\oX}{\overset{o}{X}{}}
\nc{\hl}{\overset{\leftarrow}h{}}
\nc{\hr}{\overset{\rightarrow}h{}}
\nc{\M}{{\mathcal M}}
\nc{\N}{{\mathcal N}}
\nc{\F}{{\mathcal F}}
\nc{\D}{{\mathcal D}}
\nc{\Q}{{\mathcal Q}}
\nc{\Y}{{\mathcal Y}}
\nc{\G}{{\mathcal G}}
\nc{\E}{{\mathcal E}}
\nc{\CalC}{{\mathcal C}}
\nc\Dh{\widehat{\D}}

\nc{\C}{{\mathcal C}}
\nc{\K}{{\mathcal K}}
\renewcommand{\H}{{\mathcal H}}

\nc{\T}{{\mathcal T}}
\nc{\V}{{\mathcal V}}
\renc{\P}{{\mathcal P}}
\nc{\A}{{\mathcal A}}
\nc{\B}{{\mathcal B}}
\nc{\U}{{\mathcal U}}

\nc{\Gr}{{\on{Gr}}}

\nc{\frn}{{\check{\mathfrak u}(P)}}

\nc{\fC}{\mathfrak C}
\nc{\p}{\mathfrak p}
\nc{\q}{\mathfrak q}
\nc\f{{\mathfrak f}}

\nc{\qo}{{\mathfrak q}}
\nc{\po}{{\mathfrak p}}
\nc{\s}{{\mathfrak s}}
\nc\w{\text{w}}

\renewcommand{\mod}{{\on{-mod}}}

\nc\mathi\iota
\nc\Spec{\on{Spec}}
\nc\Mod{\on{Mod}}
\nc{\tw}{\widetilde{\mathfrak t}}
\nc{\pw}{\widetilde{\mathfrak p}}
\nc{\qw}{\widetilde{\mathfrak q}}
\nc{\jw}{\widetilde j}

\nc{\grb}{\overline{\Gr}}
\nc{\I}{\mathcal I}

\nc{\lambdach}{{\check\lambda}}
\nc{\Lambdach}{{\check\Lambda}{}}
\nc{\much}{{\check\mu}}
\nc{\omegach}{{\check\omega}}
\nc{\nuch}{{\check\nu}}
\nc{\etach}{{\check\eta}}
\nc{\alphach}{{\check\alpha}}
\nc{\oblvtach}{{\check\oblvta}}
\nc{\rhoch}{{\check\rho}}
\nc{\ch}{{\check h}}

\nc{\Hb}{\overline{\H}}

\nc{\BA}{{\mathbb{A}}}
\nc{\BC}{{\mathbb{C}}}
\nc{\BF}{{\mathbb{F}}}
\nc{\BG}{{\mathbb{G}}}
\nc{\BH}{{\mathbb{H}}}
\nc{\BM}{{\mathbb{M}}}
\nc{\BO}{{\mathbb{O}}}
\nc{\BD}{{\mathbb{D}}}
\nc{\BL}{{\mathbb{L}}}
\nc{\Bl}{{\mathbb{l}}}
\nc{\BN}{{\mathbb{N}}}
\nc{\BP}{{\mathbb{P}}}
\nc{\BQ}{{\mathbb{Q}}}
\nc{\BR}{{\mathbb{R}}}
\nc{\BZ}{{\mathbb{Z}}}
\nc{\BS}{{\mathbb{S}}}

\nc{\CA}{{\mathcal{A}}}
\nc{\CB}{{\mathcal{B}}}

\nc{\CE}{{\mathcal{E}}}
\nc{\CF}{{\mathcal{F}}}
\nc{\CH}{{\mathcal{H}}}

\nc{\CL}{{\mathcal{L}}}
\nc{\CC}{{\mathcal{C}}}
\nc{\CM}{{\mathcal{M}}}
\nc{\CN}{{\mathcal{N}}}
\nc{\CK}{{\mathcal{K}}}
\nc{\CO}{{\mathcal{O}}}
\nc{\CP}{{\mathcal{P}}}
\nc{\CQ}{{\mathcal{Q}}}
\nc{\CR}{{\mathcal{R}}}
\nc{\CS}{{\mathcal{S}}}
\nc{\CT}{{\mathcal{T}}}
\nc{\CU}{{\mathcal{U}}}
\nc{\CV}{{\mathcal{V}}}
\nc{\CW}{{\mathcal{W}}}
\nc{\CX}{{\mathcal{X}}}
\nc{\CY}{{\mathcal{Y}}}
\nc{\CZ}{{\mathcal{Z}}}
\nc{\CI}{{\mathcal{I}}}
\nc{\cD}{{\mathcal{D}}}

\nc{\csM}{{\check{\mathcal A}}{}}
\nc{\oM}{{\overset{\circ}{\mathcal M}}{}}
\nc{\obM}{{\overset{\circ}{\mathbf M}}{}}
\nc{\oCA}{{\overset{\circ}{\mathcal A}}{}}
\nc{\obA}{{\overset{\circ}{\mathbf A}}{}}
\nc{\ooM}{{\overset{\circ}{M}}{}}
\nc{\osM}{{\overset{\circ}{\mathsf M}}{}}
\nc{\vM}{{\overset{\bullet}{\mathcal M}}{}}
\nc{\nM}{{\underset{\bullet}{\mathcal M}}{}}
\nc{\oD}{{\overset{\circ}{\mathcal D}}{}}
\nc{\obD}{{\overset{\circ}{\mathbf D}}{}}
\nc{\oA}{{\overset{\circ}{\mathbb A}}{}}
\nc{\op}{{\overset{\bullet}{\mathbf p}}{}}
\nc{\cp}{{\overset{\circ}{\mathbf p}}{}}
\nc{\oU}{{\overset{\bullet}{\mathcal U}}{}}
\nc{\oZ}{{\overset{\circ}{\mathcal Z}}{}}
\nc{\ofZ}{{\overset{\circ}{\mathfrak Z}}{}}
\nc{\oF}{{\overset{\circ}{\fF}}}

\nc{\fa}{{\mathfrak{a}}}
\nc{\fb}{{\mathfrak{b}}}
\nc{\fc}{{\mathfrak{c}}}
\nc{\fch}{{\mathfrak{ch}}}
\nc{\fd}{{\mathfrak{d}}}
\nc{\ff}{{\mathfrak{f}}}
\nc{\fg}{{\mathfrak{g}}}
\nc{\fgl}{{\mathfrak{gl}}}
\nc{\fh}{{\mathfrak{h}}}
\nc{\fj}{{\mathfrak{j}}}
\nc{\fl}{{\mathfrak{l}}}
\nc{\fm}{{\mathfrak{m}}}
\nc{\fn}{{\mathfrak{n}}}
\nc{\fu}{{\mathfrak{u}}}
\nc{\fp}{{\mathfrak{p}}}
\nc{\fr}{{\mathfrak{r}}}
\nc{\fs}{{\mathfrak{s}}}
\nc{\ft}{{\mathfrak{t}}}
\nc{\fz}{{\mathfrak{z}}}
\nc{\fsl}{{\mathfrak{sl}}}
\nc{\hsl}{{\widehat{\mathfrak{sl}}}}
\nc{\hgl}{{\widehat{\mathfrak{gl}}}}
\nc{\hg}{{\widehat{\mathfrak{g}}}}
\nc{\chg}{{\widehat{\mathfrak{g}}}{}^\vee}
\nc{\hn}{{\widehat{\mathfrak{n}}}}
\nc{\chn}{{\widehat{\mathfrak{n}}}{}^\vee}

\nc{\fA}{{\mathfrak{A}}}
\nc{\fB}{{\mathfrak{B}}}
\nc{\fD}{{\mathfrak{D}}}
\nc{\fE}{{\mathfrak{E}}}
\nc{\fF}{{\mathfrak{F}}}
\nc{\fG}{{\mathfrak{G}}}
\nc{\fK}{{\mathfrak{K}}}
\nc{\fL}{{\mathfrak{L}}}
\nc{\fM}{{\mathfrak{M}}}
\nc{\fN}{{\mathfrak{N}}}
\nc{\fP}{{\mathfrak{P}}}
\nc{\fU}{{\mathfrak{U}}}
\nc{\fV}{{\mathfrak{V}}}
\nc{\fZ}{{\mathfrak{Z}}}

\nc{\bb}{{\mathbf{b}}}
\nc{\bc}{{\mathbf{c}}}
\nc{\bd}{{\mathbf{d}}}
\nc{\bbf}{{\mathbf{f}}}
\nc{\be}{{\mathbf{e}}}
\nc{\bg}{{\mathbf{g}}}
\nc{\bi}{{\mathbf{i}}}
\nc{\bj}{{\mathbf{j}}}
\nc{\bn}{{\mathbf{n}}}
\nc{\bp}{{\mathbf{p}}}
\nc{\bq}{{\mathbf{q}}}
\nc{\bu}{{\mathbf{u}}}
\nc{\bv}{{\mathbf{v}}}
\nc{\bx}{{\mathbf{x}}}
\nc{\bs}{{\mathbf{s}}}
\nc{\by}{{\mathbf{y}}}
\nc{\bw}{{\mathbf{w}}}
\nc{\bA}{{\mathbf{A}}}
\nc{\bK}{{\mathbf{K}}}
\nc{\bB}{{\mathbf{B}}}
\nc{\bC}{{\mathbf{C}}}
\nc{\bG}{{\mathbf{G}}}
\nc{\bD}{{\mathbf{D}}}
\nc{\bH}{{\mathbf{He}}}
\nc{\bM}{{\mathbf{M}}}
\nc{\bN}{{\mathbf{N}}}
\nc{\bO}{{\mathbf{O}}}
\nc{\bV}{{\mathbf{V}}}
\nc{\bW}{{\mathbf{Wh}}}
\nc{\bX}{{\mathbf{X}}}
\nc{\bZ}{{\mathbf{Z}}}
\nc{\bS}{{\mathbf{S}}}
\nc{\bT}{{\mathbf{T}}}

\nc{\sA}{{\mathsf{A}}}
\nc{\sB}{{\mathsf{B}}}
\nc{\sC}{{\mathsf{C}}}
\nc{\sD}{{\mathsf{D}}}
\nc{\sF}{{\mathsf{F}}}
\nc{\sG}{{\mathsf{G}}}
\nc{\sK}{{\mathsf{K}}}
\nc{\sM}{{\mathsf{M}}}
\nc{\sO}{{\mathsf{O}}}
\nc{\sU}{{\mathsf{U}}}
\nc{\sW}{{\mathsf{W}}}
\nc{\sQ}{{\mathsf{Q}}}
\nc{\sP}{{\mathsf{P}}}
\nc{\sZ}{{\mathsf{Z}}}
\nc{\sfp}{{\mathsf{p}}}
\nc{\sfq}{{\mathsf{q}}}
\nc{\sr}{{\mathsf{r}}}
\nc{\sk}{{\mathsf{k}}}
\nc{\su}{{\mathsf{u}}}
\nc{\sv}{{\mathsf{v}}}
\nc{\sg}{{\mathsf{g}}}
\nc{\sff}{{\mathsf{f}}}
\nc{\sfb}{{\mathsf{b}}}
\nc{\sfc}{{\mathsf{c}}}
\nc{\sd}{{\mathsf{d}}}

\nc{\BK}{{\bar{K}}}

\nc{\tA}{{\widetilde{\mathbf{A}}}}
\nc{\tB}{{\widetilde{\mathcal{B}}}}
\nc{\tg}{{\widetilde{\mathfrak{g}}}}
\nc{\tG}{{\widetilde{G}}}
\nc{\TM}{{\widetilde{\mathbb{M}}}{}}
\nc{\tO}{{\widetilde{\mathsf{O}}}{}}
\nc{\tU}{{\widetilde{\mathfrak{U}}}{}}
\nc{\TZ}{{\tilde{Z}}}
\nc{\tx}{{\tilde{x}}}
\nc{\tbv}{{\tilde{\bv}}}
\nc{\tfP}{{\widetilde{\mathfrak{P}}}{}}
\nc{\tz}{{\tilde{\zeta}}}
\nc{\tmu}{{\tilde{\mu}}}

\nc{\urho}{\underline{\rho}}
\nc{\uB}{\underline{B}}
\nc{\uC}{{\underline{\mathbb{C}}}}
\nc{\ui}{\underline{i}}
\nc{\uj}{\underline{j}}
\nc{\ofP}{{\overline{\mathfrak{P}}}}
\nc{\oB}{{\overline{\mathcal{B}}}}
\nc{\og}{{\overline{\mathfrak{g}}}}
\nc{\oI}{{\overline{I}}}

\nc{\eps}{\varepsilon}
\nc{\hrho}{{\hat{\rho}}}

\nc{\one}{{\mathbf{1}}}
\nc{\two}{{\mathbf{t}}}

\nc{\Rep}{{\mathop{\operatorname{\rm Rep}}}}
\nc{\Tot}{{\mathop{\operatorname{\rm Tot}}}}
\nc{\Ker}{{\mathop{\operatorname{\rm Ker}}}}
\nc{\Hilb}{{\mathop{\operatorname{\rm Hilb}}}}
\nc{\End}{{\mathop{\operatorname{\rm End}}}}
\nc{\Ext}{{\mathop{\operatorname{\rm Ext}}}}
\nc{\CHom}{{\mathop{\operatorname{{\mathcal{H}}\it om}}}}
\nc{\GL}{{\mathop{\operatorname{\rm GL}}}}
\nc{\gr}{{\mathop{\operatorname{\rm gr}}}}
\nc{\Id}{{\mathop{\operatorname{\rm Id}}}}
\nc{\de}{{\mathop{\operatorname{\rm def}}}}
\nc{\length}{{\mathop{\operatorname{\rm length}}}}
\nc{\supp}{{\mathop{\operatorname{\rm supp}}}}

\nc{\Cliff}{{\mathsf{Cliff}}}
\nc{\Fl}{\on{Fl}}
\nc{\Fib}{{\mathsf{Fib}}}
\nc{\Coh}{{\on{Coh}}}
\nc{\QCoh}{{\on{QCoh}}}
\nc{\IndCoh}{{\on{IndCoh}}}
\nc{\FCoh}{{\mathsf{FCoh}}}

\nc{\reg}{{\text{\rm reg}}}

\nc{\cplus}{{\mathbf{C}_+}}
\nc{\cminus}{{\mathbf{C}_-}}
\nc{\cthree}{{\mathbf{C}_*}}
\nc{\Qbar}{{\bar{Q}}}
\nc\Eis{{\on{Eis}}}
\nc\Eisb{\ol\Eis{}}
\nc\Eisr{\on{Eis}^{rat}{}}
\nc\wh{\widehat}
\nc{\Def}{\on{Def_{\check{\fb}}(E)}}
\nc{\barZ}{\overline{Z}{}}
\nc{\barbarZ}{\overline{\barZ}{}}
\nc{\barpi}{\overline\pi}
\nc{\barbarpi}{\overline\barpi}
\nc{\barpip}{\overline\pi{}^+}
\nc{\barpim}{\overline\pi{}^-}

\nc{\fq}{\mathfrak q}

\nc{\fqb}{\ol{\fq}{}}
\nc{\fpb}{\ol{\fp}{}}
\nc{\fpr}{{\fp^{rat}}{}}
\nc{\fqr}{{\fq^{rat}}{}}

\nc{\hattimes}{\wh\otimes}

\nc{\bh}{{{\mathbf h}}}
\nc{\bk}{{{\mathbf k}}}
\nc{\bOmega}{{\overline{\Omega(\check \fn)}}}

\nc{\seq}[1]{\stackrel{#1}{\sim}}

\nc{\cT}{{\check{T}}}
\nc{\cG}{{\check{G}}}
\nc{\cM}{{\check{M}}}
\nc{\cB}{{\check{B}}}
\nc{\cP}{{\check{P}}}

\nc{\ct}{{\check{\mathfrak t}}}
\nc{\cg}{{\check{\fg}}}
\nc{\cb}{{\check{\fb}}}
\nc{\cn}{{\check{\fn}}}

\nc{\cLambda}{{\check\Lambda}}

\nc{\cla}{{\check\lambda}}
\nc{\cmu}{{\check\mu}}
\nc{\cnu}{{\check\nu}}
\nc{\ceta}{{\check\eta}}

\nc{\DefbE}{{\on{Def}_{\cB}(E_\cT)}}

\nc{\imathb}{{\ol{\imath}}}
\nc{\rlr}{\overset{\longrightarrow}{\underset{\longrightarrow}\longleftarrow}}

\nc{\oBun}{\overset{\circ}\Bun}
\nc{\LocSys}{\on{LocSys}}
\nc{\BunBbb}{\ol{\ol{Bun}}_B}
\nc{\BunBr}{\Bun_B^{rat}}
\nc{\BunBrp}{\Bun_B^{rat,polar}}
\nc{\BunTrp}{\Bun_T(X)^{rat,polar}}
\nc{\BunNr}{\Bun_N^{rat}}
\nc{\BunNre}{\Bun_N^{enh,rat}}
\nc{\BunTr}{\Bun_T(X)^{rat}}
\nc{\Vect}{\on{Vect}}
\nc{\Whit}{\on{Whit}}
\nc{\CTb}{\ol{\on{CT}}}
\nc{\Ran}{\on{Ran}}
\nc{\CTr}{\on{CT}^{rat}{}}
\nc\jmathr{\jmath^{rat}{}}
\nc{\ux}{\underline{x}}
\nc{\clambda}{{\check\lambda}}
\nc{\calpha}{{\check\alpha}}
\nc{\ind}{{\mathbf{ind}}}
\nc{\oblv}{{\mathbf{oblv}}}
\nc{\coeff}{\on{W-coeff}}
\nc{\Poinc}{\on{Poinc}}
\nc{\Dmod}{\on{D}}
\nc{\dr}{\on{dR}}
\nc{\oCZ}{\overset{\circ}\CZ}
\nc{\KL}{\on{KL}}
\nc{\triv}{{\mathbf{triv}}}

\nc{\dgSch}{\on{DGSch}}
\nc{\Sch}{\on{Sch}}
\nc{\affdgSch}{\on{Sch}^{\on{aff}}}
\nc{\affSch}{\on{Sch}^{\on{aff}}}
\nc{\Sing}{\on{Sing}}
\nc{\inftygroup}{\infty\on{-Grpd}}
\renc{\dr}{{\on{dr}}}
\nc\Maps{\on{Maps}}
\nc\bMaps{\mathbf{Maps}}
\nc{\ul}{\underline}
\nc{\bNP}{\mathbf{N(P)}}
\nc{\ofc}{\overset{\circ}\fch}
\nc{\ppart}{(\!(t)\!)}
\nc{\qqart}{[\![t]\!]}
\nc{\crit}{\on{crit}}
\nc{\bDelta}{\mathbf{\Delta}}
\nc{\genB}{{\overset{\on{gen}}\to B}}
\nc{\genP}{{\underset{\on{gen}}\longrightarrow P}}
\nc{\genN}{{\underset{\on{gen}}\longrightarrow N}}
\nc{\sotimes}{\overset{!}\otimes}
\nc{\mmod}{{\on{-}\mathbf{mod}}}
\nc{\Spc}{\on{Spc}}
\nc{\bLoc}{{\mathbf {Loc}}}
\nc{\bGamma}{{\mathbf \Gamma}}

\begin{document}

\title[Recent progress in geometric Langlands theory]{Recent progress in geometric Langlands theory}

\author{Dennis Gaitsgory} 

\begin{abstract}
The is the English version of the text of the talk at S\'eminaire Bourbaki on February 16, 2016.
\end{abstract} 

\date{\today}

\maketitle

\section{Introduction}

Throughout the talk we fix $X$ to be a smooth connected complete curve and $G$ a reductive group
over a ground field $k$.

\medskip

When discussing connections with the classical (function-theoretic) Langlands theory, 
we will assume that $k=\BF_q$. When talking about the \emph{categorical} geometric
Langlands theory, we will take $k$ to be characteristic zero. 

\ssec{Some history}

What is nowadays knows as the \emph{geometric Langlands theory} originated from the ideas of four people: 
A.~Beilinson, P.~Deligne, V.~Drinfeld and G.~Laumon. 

\sssec{}

The first input was Deligne's observation that one can prove the existence of the grossen-character corresponding to a
unramified character of the Galois of a function field using algebro-geometric considerations. The idea is the following.

\medskip

An (unramified) grossen-character can be thought of as a function on the set (rather, groupoid) of $\BF_q$-points of the Picard stack $\on{Pic}(X)$ of $X$.
We will do the construction in two steps. First, starting from an unramified Galois character $\sigma$, we will construct an 
$\ell$-adic sheaf on $\on{Pic}(X)$, denoted $\CF_\sigma$. Then the sought-for grossen-character will be obtained from $\CF_\sigma$
by Grothedieck's sheaves-functions correspondence, i.e., by taking traces of the Frobenius. 

\medskip

The construction of $\CF_\sigma$ is geometric. Namely, we interpret $\sigma$ as a 1-dimensional $\ell$-adic local system on 
$X$, denoted $E_\sigma$. To $E_\sigma$ and $d\geq 0$, we attach the \emph{symmetric power}, denoted $E_\sigma^{(d)}$,
which is a 1-dimensional local system on the scheme $X^{(d)}$ parameterizing effective divisors on $X$ of degree $d$. 
(The sheaf $E_\sigma^{(d)}$ is a natural thing to do from the number-theoretic point of view: the function attached to it is the function
corresponding to $\sigma$ on the set of effective divisors.)

\medskip

Now, we consider the Abel-Jacobi map $$X^{(d)}\to \on{Pic}(X)$$
and the task is to show that there exists an $\ell$-adic sheaf $\CF_\sigma$ on $\on{Pic}(X)$ that pulls back to $E_\sigma^{(d)}$
for each $d$. It is easy to that it is enough to prove the existence of $\CF_\sigma$ over the connected components
$\on{Pic}^d(X)$ of $\on{Pic}(X)$ for $d$ large (i.e., $d\geq d_0$ for some fixed $d_0$).  

\medskip

The punchline is that for $d>2g-2$
(here $g$ is the genus of $X$),
the map $X^{(d)}\to \on{Pic}^d(X)$
is a smooth fibration with simply-connected fibers, which guarantees the existence (and uniqueness) of the descent of
$E_\sigma^{(d)}$ to the sought-for $\ell$-adic sheaf $\CF_\sigma^d$ on $\on{Pic}^d(X)$.

\sssec{}

Then came Drinfeld's ground-breaking paper \cite{Dr}. In a sense it was the extension of Deligne's construction to the vastly more complicated
case, when instead of grossen-characters we consider unramified automorphic functions for the group $GL_2$.  Here again, we interpret the
(unramified) automorphic space as the groupoid of $\BF_q$-points of the moduli space $\Bun_2:=\Bun_{GL_2}$ classifying
rank-2 vector bundles on $X$. Drinfeld's idea is to attach to a 2-dimensional Galois representation $\sigma$ an $\ell$-adic sheaf
(by which we actually mean an object of the corresponding derived category) $\CF_\sigma$ on $\Bun_2$, and then obtain the sought-for function
by taking the traces of the Frobenius.

\medskip

The main difference from the commutative case, considered by Deligne (which corresponds to the case of the group $\BG_m=GL_1$), 
is that the construction of $\CF_\sigma$ starting from $\sigma$ is much more involved. The intermediate player, i.e., $E_\sigma^{(d)}$,
is now interpreted as the $\ell$-adic sheaf that records the \emph{Whittaker} (a.k.a., Fourier) 
coefficients of $\CF_\sigma$. So, our task is to reconstruct an automorphic object from its Fourier coefficients. This is again done via
appealing to geometry--ultimately the simply-connectedness of fibers of some map. 

\sssec{}

After Drinfeld's paper came one by Laumon, \cite{Lau1}, which gave a conjectural extension of Drinfeld's construction from $GL_2$
to $GL_n$. To the best of our knowledge, the title of Laumon's paper was the first place where the combination of words 
`geometric Langlands' appeared. 

\medskip

While the stated goal of Drinfeld's paper was to construct an automorphic \emph{function}, Laumon's paper had the effect of
shifting the goal: people became interested in automorphic \emph{sheaves} ($\ell$-adic sheaves on $\Bun_n(X)$) for their own sake.

\medskip

Following the appearance of Laumon's paper, it became clear that one should also try to attack $\Bun_G(X)$ for an
arbitrary reductive $G$, even though it was not clear how to do this (because the Whittaker model does not work as nicely
outside the case of $G=GL_n$). 

\sssec{}

The next paradigm shift came in the work of Beilinson and Drinfeld, \cite{BD}. They considered the same $\Bun_G(X)$, but now
over a ground field $k$ of characteristic zero, and instead of $\ell$-adic sheaves, they proposed to consider D-modules.

\medskip

In this case, a new method for constructing objects becomes available: by generators and relations.  A fancy version of `generators and
relations' principle--the localization functor, pioneered in \cite{BB}, lies in the core of the manuscript \cite{BD}, which produces
automorphic D-modules using representations of the Kac-Moody Lie algebra (thought of as the Lie algebra of infinitesimal
symmetries of a $G$-bundle on the formal punctured disk). 

\sssec{}

In an independent development, in \cite{Lau2}, Laumon showed that if we take $G$ to be a torus $T$, a generalized version
of the Fourier-Mukai transform identifies the (derived) category of D-modules on the stack $\Bun_T(X)$ with the (derived) category
of quasi-coherent sheaves on the stack $\on{LocSys}_{\check{T}}(X)$ of \emph{de Rham local systems} on $X$ with respect to the
Langlands dual torus $\cT$. 

\medskip

I.e., Laumon's paper extends the poinwtise Langlands correspondence (i.e., construction of $\CF_\sigma$ corresponding to a fixed
local system $\sigma$) to a statement about the universal family of local systems.

\sssec{}  \label{sss:best hope intro}

Finally, combining Laumon's equivalence for the torus, and accumulated evidence for the general $G$, Beilinson and Drinfeld came
up with the idea of \emph{categorical geometric Langlands equivalence}. 

\medskip

In its crude form, this should be an equivalence between the 
(derived) category $\Dmod(\Bun_G(X))$ of D-modules on the stack $\Bun_G(X)$ and the (derived) category $\QCoh(\LocSys_\cG(X))$
of quasi-coherent sheaves on the stack 
$\on{LocSys}_{\check{G}}(X)$. Such as equivalence is what Beilinson and Drinfeld called the \emph{best hope}, but they 
never stated it explicitly, because it is (and was) known that it cannot hold \emph{as-is} beyond the case of a torus (the reason reason 
for this will be indicated in \secref{sss:arth}). 

\ssec{What do we mean by `geometric Langlands' nowadays?}

There are several meta-problems that comprise what one can call the \emph{geometric Langlands theory}; we shall list some
of them below; the order in which they will appear reflects (our perception of) the historical development (and the increasing level of technical complexity) rather
than how the complete picture should ultimately look like (e.g., we do think that the quantum case is more fundamental than the usual one). 

\medskip

We will only consider the \emph{categorical} geometric Langlands theory; in particular we will assume that the ground field $k$ is of characteristic 
zero, and on the automorphic side we will work with D-modules rather than $\ell$-adic sheaves. 

\medskip

We should remark that whatever conjectures and meta-conjectures we mention below, they are all theorems when the group
$G$ is a torus, thanks to the various generalizations of the Deligne-Fourier-Mukai-Laumon transform. 

\sssec{}

First, we have the \emph{categorical\footnote{From now on, we will drop the adjective `categorical', because everything will be categorical. }
global unramified} geometric Langlands. This is an attempt to formulate and prove
a version of the \emph{best hope} by Beilinson and Drinfeld, mentioned above. I.e., we want a category that it a close
cousin (or identical twin) of $\Dmod(\Bun_G(X))$ to be equivalent to a category that is a close cousin of $\QCoh(\on{LocSys}_{\check{G}}(X))$.

\medskip

This is the aspect of the geometric Langlands theory that has been developed the most. It will be discussed in Sects. \ref{s:Hecke} and \ref{s:global}.

\sssec{}

Next there is the \emph{local ramified} geometric Langlands theory. Unlike the global case, in the local version we are interested in an equivalence of
\emph{2-categories} (rather than 1-categories, i.e., just categories). 
For a long time it was not even clear how to formulate our wish (specifically, what 2-category to consider on the Galois side). 
However, recently, a breakthrough has been achieved in the work of S.~Raskin, \cite{Ras}. We will discuss this in \secref{s:loc}. 

\medskip

We should also mention that the tamely ramified case of the local ramified geometric Langlands had been settled by R.~Bezrukavnikov in \cite{Bez}
even before the general program was formulated. 

\sssec{}

Next, there is the global ramified Langlands theory. Its tamely ramified case has not been explicitly studied in detail, but the current state of
knowledge should allow to bring it to the same status as the unramified case.

\medskip

The general ramified case is wide-open, and there are formidable technical difficulties that one needs to surmount in order to start investigating it.
One of the difficulties is that we do not know whether the category of D-modules on the automorphic side, i.e., the (derived) category of D-modules
on the moduli space $\Bun_G(X)^{k\cdot x}$
of $G$-bundles on $X$ \emph{equipped with structure of level $k\geq 1$ at a point} $x$ is 
compactly generated\footnote{If one surveys the literature, in most of the statements that involve an equivalences of two triangulated/DG 
categories, the categories of question are compactly generated. The reason is that we do not know very well how to compute things outside the 
compactly generated case.}.

\sssec{}

Finally, all of the above three aspects: unramified global, ramified local and ramified global admit \emph{quantum} versions. The quantum parameter 
in the \emph{quantum geometric Langlands theory} is a non-degenerate $W$-invariant symmetric bilinear form on the Cartan Lie algebra
$\fh$ of $\fg$, whose inverse is a similar kind of datum for $\cg$. 

\medskip

At the risk of making a controversial statement, the author/speaker has to admit that  
he came to regard the quantum version as the ultimate reason of `why something like the Langlands theory takes place', with the usual geometric 
Langlands being its degeneration (letting the quantum parameter tend to zero), and the classical (i.e., function-theoretic) Langlands theory as some 
sort of residual phenomenon. 

\medskip 

We will not say anything about the quantum case in this talk, but rather refer the reader to \cite{Ga2}, where the dream of \emph{quantum geometric Langlands} 
is discussed. Here we will only mention the following two facts. 

\medskip

One is that in the quantum theory restores the symmetry between $G$ and $\cG$: we no longer have the 
Galois side, but rather both sides are automorphic, but twisted by the quantum parameter. 

\medskip

The other is that the guiding principle of the quantum theory is 
that `Whittaker is dual to Kac-Moody', which is striking because `Whittaker' has a classical (i.e., number-theoretic) meaning, while `Kac-Moody' does not.

\ssec{Terminology and notation}

\sssec{}

The global unramified geometric Langlands conjecture can be formulated as an equivalence of (triangulated) categories.
But if one wants to dig a tiny bit deeper into attempts of its proof, one needs to work with $\infty$-categories--in this
case with ($k$-linear) DG categories. We refer the reader to \cite[Sect. 10]{GR1} for the definition of the latter. 

\medskip

For the reader not familiar with $\infty$-categories, we recommend the following approach. On the first pass pretend that there no
difference between $\infty$-categories and ordinary categories. On the second pass pretend that you already know what $\infty$-categories are
and stay tuned for the language used when working with them (a survey of the syntax of $\infty$-categories can be found in 
\cite[Sect. 1]{Lu}, or from a somewhat different perspective, in \cite[Sect. 1]{GR1}). On the third pass...learn the theory properly!

\sssec{}

Another piece of  `bad news' is that when working on the Galois side of the geometric Langlands theory, we cannot stay
within the realm of classical algebraic geometry, and one needs to plunge oneself into the world of DAG--derived algebraic
geometry. For example, the stack $\LocSys_\cG(X)$ has a non-trivial derived structure for $G=T$ being a torus. 
The reader is referred to \cite{GR2} for an introduction to DAG. 

\medskip

In what follows, when we say `scheme' or `algebraic stack', we will tacitly mean the corresponding derived notions.

\sssec{}

To a scheme or algebraic stack $Y$ one attaches the DG category $\QCoh(Y)$; we will somewhat abusively refer to it as the 
\emph{(derived) category of quasi-coherent sheaves on}\footnote{When $Y$ is a \emph{classical} (as opposed to derived) scheme
or a sufficiently nice algebraic stack, this category \emph{is} the derived category of its heart with respect to a naturally defined t-structure.} 
$Y$. We refer the reader to \cite{GR3} for the definition. 
We note, however, that the definition of $\QCoh(Y)$ is much more general: it makes sense for $Y$ which is an arbitrary \emph{prestack}\footnote{
But in this more general setting, $\QCoh(Y)$ is not at all the \emph{derived} category of any abelian category, even if $Y$ itself is a classical
prestack.}. 

\sssec{}

If $Y$ is a scheme/algebraic stack/prestack locally of finite type over $k$, one can attach to it the DG category $\on{D}(Y)$;
we will also somewhat abusively refer to it as the \emph{(derived) category of D-modules on}\footnote{If $Y$ is a scheme,
this \emph{is} the derived category of its heart, but this is no longer true for algebraic stacks.} $Y$. We refer the reader
to \cite{GR4} for the definition. 

\medskip

The good news is that when discussing $\on{D}(Y)$, the derived structure on $Y$ plays no role. So, on the automorphic
side of the geometric Langlands theory we can stay within the realm of classical algebraic geometry. 

\sssec{}

Whenever we talk about functors between (derived) categories of sheaves/D-modules on various spaces (!- and *- direct and inverse images),
we will always mean the corresponding derived functors. I.e., abelian categories will not appear unless explicitly stated
otherwise. 

\sssec{}

For the motivational parts of the talk (i.e., analogies with the function-theoretic situation), we will assume the reader's familiarity
with the basics of algebraic number theory (ad\`eles, ramification, Frobenius elements, etc.) 

\ssec{Acknowledgements}

The author/speaker wishes to thank D.~Arinkin, A.~Beilinson, J.~Bernstein, V.~Drinfeld, E.~Frenkel, D.~Kazhdan, 
S.~Raskin and E.~Witten, discussions with whom has informed his perception of the Langlands theory. 

\section{Hecke action}  \label{s:Hecke}

The classical Langlands correspondence, and historically also the geometric one, were characterized by relating the
spectrum of the action of the Hecke operators (resp., functors) on the automorphic side to a Galois datum. 
We begin by discussing this aspect of the theory. 

\ssec{Hecke action on automorphic functions}

Let $\CK$ be the field of rational functions on our curve $X$; let $\BA$ be the ring of ad\`eles, and $\BO\subset \BA$
the subring of integral ad\`eles. For a place $x\in X$, we let $\CO_x\subset \CK_x$ denote the corresponding local
ring and local field, respectively. 

\sssec{}

The automorphic space is by definition the quotient $G(\CA)/G(\CK)$. It is acted on by left translations by the group
$G(\CA)$. The unramified automorphic space is the set (but, properly speaking, \emph{groupoid})
$$G(\BO)\backslash G(\BA)/G(\CK).$$

\medskip

Our object of study is the space $\on{Autom}(X)$ of unramified $\BQ_\ell$-valued automorphic functions,
i.e., functions on $G(\BO)\backslash G(\BA)/G(\CK)$, or, which is the same, the space of
$G(\BO)$-invariant functions on $G(\CA)/G(\CK)$. 

\sssec{}

Since the subgroup $G(\BO)\subset G(\BA)$ is not normal, we do not have an action of $G(\BA)$ on $\on{Autom}(X)$. Instead, 
the action of $G(\BA)$ on $G(\BA)/G(\CK)$ induces an action on $\on{Autom}(X)$ of the \emph{spherical Hecke algebra}
$\BH(G)_X$. By definition, as a vector space, $\BH(G)_X$ consists of compactly supported $G(\BO)$-biinvariant functions on $G(\BA)$,
and it is endowed with a structure of associative algebra via the operation of \emph{convolution}.

\medskip

The datum of the action of $\BH(G)_X$ is equivalent to that of a family of pairwise commuting actions of the local Hecke
algebras $\BH(G)_x$ for every place $x$ of $X$, where each $\BH(G)_x$ is the algebra (with respect to convolution)
of $G(\CO_x)$-biinvariant compactlt supported functions on $G(\CK_x)$.

\medskip

Our interest is to find the spectrum of $\BH(G)_X$ (i.e., the joint spectrum of the algebras $\BH(G)_x$) acting on $\on{Autom}(X)$. 

\sssec{}  \label{sss:Satake}

Fix $x\in X$. The first basic fact about the associative algebra $\BH(G)_x$ is that it is actually commutative. But, in addition to this, we can
actually describe it very explicitly. 

\medskip

Namely, the classical Satake isomorphism says that $\BH(G)_x$ identifies canonically with the algebra of
ad-invariant regular functions on the algebraic group $\cG$, where $\cG$ is the Langlands dual of $G$, thought of as an algebraic
group over $\BQ_\ell$.

\sssec{}

The Satake isomorphism allows us to give a formulation of \emph{Langlands correspondence}.

\medskip

Namely, given a unramified representation $\sigma$ of the Galois group of $\CK$ into $\cG$ (defined up to conjuagation)
and a common eigenvector $f\in \on{Autom}(X)$ of the algebras $\BH(G)_x$, we shall say that $f$ corresponds to 
$\sigma$, if for every $x\in X$, the character by which $\BH(G)_x$ acts on $f$ is given, in terms of the Satake isomorphism,
by \emph{evaluation} of functions on $\cG$ on the conjugacy class of the image of $\on{Frob}_x$ under the map $\sigma$.
Here $\on{Frob}_x$ is the Frobenius at $x\in X$, which is a well-defined conjugacy class in the unramified quotient of
the Galois group of $\CK$. 

\sssec{}  \label{sss:Vinc}

The recent result of V.~Lafforgue (known as the Automorphic $\Rightarrow$ Galois direction of Langlands correspondence, see \cite{VLaf})
says that for every eigenvector $f$ (assumed cuspidal) \emph{there exists} a $\sigma$ to which this $f$ corresponds. In this case of
$G=GL_n$ the existence statement can be strengthened to one about uniqueness and surjectivity, due to the work of L.~Lafforgue \cite{LLaf}. 

\ssec{Geometric Satake and Hecke eigensheaves}

We now pass to considering the Hecke action in the geometric context, i.e., when instead of functions on the automorphic
space we consider the appropriately defined derived category $\on{D}(\Bun_G(X))$ of $\ell$-adic sheaves/D-modules on the automorphic
stack $\Bun_G(X)$. 

\sssec{}

The initial observation is the \emph{geometric Satake} equivalence of Lusztig-Drinfeld-Ginzburg-Mirkovi\'c-Vilonen (historical order)
that says that for every point $x\in X$ the monoidal category $\Rep(\cG)$ of algebraic representations of $\cG$ acts on the category $\on{D}(\Bun_G(X))$. 

\medskip

Thus, we obtain the Hecke functors
$$\on{H}_{V,x}:\on{D}(\Bun_G(X))\to \on{D}(\Bun_G(X)),\quad x\in X,\,\,V\in \Rep(\cG).$$

\medskip

This is the geometric replacement of the $\BH(G)_x$-action on $\on{Autom}(X)$ combined with the Satake isomorphism of \secref{sss:Satake}. 

\medskip

However, in geometry one can do much more: one can make the point $x$ move along $X$. Thus, for every $V\in \Rep(\cG)$ we obtain the Hecke
functor
$$\on{H}_{V}:\on{D}(\Bun_G(X))\to \on{D}(\Bun_G(X)\times X).$$

\sssec{}

But in fact, one can do even more than that. Let is take a pair of objects $V_1,V_2\in \Rep(\cG)$. To them we can canonically attach the functor
$$\on{H}_{V_1,V_2}:\on{D}(\Bun_G(X))\to \on{D}(\Bun_G(X)\times X\times X),$$
which is the composition
$$\on{D}(\Bun_G(X))\overset{\on{H}_{V_1}}\longrightarrow \on{D}(\Bun_G(X)\times X) \overset{\on{H}_{V_2}}\longrightarrow \on{D}(\Bun_G(X)\times X\times X)$$
and \emph{also} the composition
$$\on{D}(\Bun_G(X))\overset{\on{H}_{V_2}}\longrightarrow \on{D}(\Bun_G(X)\times X) \overset{\on{H}_{V_1}}\longrightarrow \on{D}(\Bun_G(X)\times X\times X),$$
up to the permutation of the factors in $X\times X$.

\medskip

A key property of this functor is that it composition with the restriction functor
$$(\on{id}_{\Bun_G(X)}\times \on{diag}_X)^!:\on{D}(\Bun_G(X)\times X\times X)\to \on{D}(\Bun_G(X)\times X)$$
identifies with the functor $\on{H}_{V_1\otimes V_2}$. 

\sssec{}  \label{sss:VI}

Let us now say the same but slightly more generally and abstractly. Let $I$ be a non-empty finite set, and let $V_I$ be an $I$-tuple of
objects of $\Rep(\cG)$
$$(i\in I)\rightsquigarrow V_i\in \Rep(\cG).$$

To this datum we attach a functor
$$\on{H}_{V_I}:\on{D}(\Bun_G(X))\to \on{D}(\Bun_G(X)\times X^I).$$

When $I$ is a singleton, we recover the functor $\on{H}_V$. 

\medskip

The assignment $I\mapsto \on{H}_{V_I}$ is compatible with 
the operation of disjoint union of finite sets: for $I=I_1\sqcup I_2$, the functor $\on{H}_{V_I}$ identifies with
$$\on{D}(\Bun_G(X))\overset{\on{H}_{V_{I_1}}}\longrightarrow \on{D}(\Bun_G(X)\times X^{I_1}) 
\overset{\on{H}_{V_{I_2}}}\longrightarrow \on{D}(\Bun_G(X)\times X^{I_1}\times X^{I_2})\simeq \on{D}(\Bun_G(X)\times X^I).$$

\sssec{}

Let us now be given a surjection of finite sets $\phi:I\twoheadrightarrow J$. Given $V_I:I\to \Rep(\cG)$ we can create 
$\phi(V_I)=:V_J$ by
\begin{equation} \label{e:VJ}
V_j=\underset{i\in \phi^{-1}(j)}\bigotimes\, V_i.
\end{equation} 

Let $\on{diag}_\phi$ denote the map $X^J\to X^I$, corresponding to $\phi$. Then the composition
$$\on{D}(\Bun_G(X))\overset{\on{H}_{V_I}} \longrightarrow\on{D}(\Bun_G(X)\times X^I)\overset{(\on{id}_{\Bun_G(X)}\times \on{diag}_\phi)^!}
\longrightarrow \on{D}(\Bun_G(X)\times X^J)$$
identifies with $\on{H}_{V_J}$. 

\sssec{}

We will now perform one more manipulation. Let $\CM_I$ be an object of $\on{D}(X^I)$. We define the \emph{endo-functor}
$\on{H}_{V_I,\CM_I}$ of $\on{D}(\Bun_G(X))$ to be the composition of $\on{H}_{V_I}$, followed by the functor 
$$\on{D}(\Bun_G(X)\times X^I)\to \on{D}(\Bun_G(X)), \quad \CF\mapsto (\on{pr}_{\Bun_G(X)})_!(\CF\sotimes (\on{pr}_{X^I})^!(\CM_I)).$$
Here $\on{pr}_{\Bun_G(X)}$ and $\on{pr}_{X^I}$ denote the two projections
$$\Bun_G(X) \leftarrow \Bun_G(X)\times X^I\rightarrow X^I,$$
and $\sotimes$ is the !-tensor product of sheaves/D-modules (the !-pullback of the external tensor product by the diagonal morphism). 

\sssec{}

For $I=I_1\sqcup I_2$ and $\CM=\CM_1\boxtimes \CM_2\in \on{D}(X^I)\simeq \on{D}(X^{I_1}\times X^{I_2})$, we have
$$\on{H}_{V_I,\CM_I}\simeq \on{H}_{V_{I_1},\CM_{I_1}}\circ \on{H}_{V_{I_2},\CM_{I_2}}\simeq 
\on{H}_{V_{I_2},\CM_{I_2}}\circ \on{H}_{V_{I_1},\CM_{I_1}}.$$

\medskip

For $\phi:I\twoheadrightarrow J$ and $\CM_J\in \on{D}(X^J)$ we have
$$\on{H}_{V_I,(\on{diag}_\phi)_!(\CM_J)}\simeq \on{H}_{V_J,\CM_J}.$$

\sssec{}

Now, the collection of $(I,V_I,\CM_I)$ can be glued to a category\footnote{As always, `category' means `DG category'.},
by imposing the following family of relations\footnote{Formally,
we take the co-end in $\on{DGCat}$ of the following functors from the category of non-empty finite sets
and surjections: one takes $I$ to $\Rep(\cG)^{\otimes I}$ and another to $\on{D}(X^I)$.} 

$$(I,V_I,\CM_I)\simeq (J,V_J,\CM_J)$$
each time we have 
$$\phi:I\twoheadrightarrow J, \quad \CM_I=(\on{diag}_\phi)_!(\CM_J), \quad V_J=\phi(V_I).$$

We denote this category by $\Rep(\cG,\Ran(X))$ (see \cite[Sect. 4.2]{Ga3} for another approach 
to defining $\Rep(\cG,\Ran(X))$). 

\medskip

The operation of disjoint union of finite sets induces on $\Rep(\cG,\Ran(X))$ a structure of \emph{non-unital} (symmetric)
monoidal category. 

\sssec{}  \label{sss:Hecke act ult}

The upshot of all the preceding discussion is that, ultimately, the geometric Hecke action amounts to the action of
the monoidal category $\Rep(\cG,\Ran(X))$ on  $\on{D}(\Bun_G(X))$. 

\ssec{The spectral decomposition}

In this subsection we will assume that $k$ is of characteristic zero, and we will take $\on{D}(-)$ to mean the derived category 
of D-modules.

\medskip

We consider the stack $\LocSys_\cG(X)$ of de Rham $\cG$-local systems on $X$. We will explain the counterpart within 
the geometric Langlands theory of V.~Lafforgue's theorem mentioned in \secref{sss:Vinc}, which heuristically means that the
spectrum of the Hecke functors on $\on{D}(\Bun_G(X))$ is contained in $\LocSys_\cG(X)$.

\sssec{}

First we fix a point $x\in X$ and an object $V\in \Rep(\cG)$. To this data we associate a coherent sheaf (in fact, a vector bundle)
on $\LocSys_\cG(X)$, denoted $\on{Ev}_{V,x}$. 

\medskip

Namely, the fiber of $\on{Ev}_{V,x}$ at a point $\sigma\in \LocSys_\cG(X)$ is $(V^\sigma)_x$, where $V^\sigma$ is the local system
(=lisse D-module) associated to the $\cG$-representation $V$ and the $\cG$-local system $\sigma$, and $(-)_x$ denotes taking
the !-fiber at $x$. 

\medskip

More generally, given a finite set $I$ and $V_I$ as in \secref{sss:VI}, we can associate to this data an object $\on{Ev}_{V_I}$, which
is a quasi-coherent sheaf on $\LocSys_\cG(X)\times X^I$,  \emph{equipped with a connection along} $X^I$.

\medskip

Hence, given in addition $\CM_I\in \on{D}(X^I)$, we can produce 
\begin{equation} \label{e:Ev}
\on{Ev}_{V_I,\CM_I}\in \QCoh(\LocSys_\cG(X))
\end{equation} 
by taking the de Rham direct image of $\on{Ev}_{V_I}\otimes (\on{pr}_{X^I})^!(\CM_I)$ along the projection
$$\on{pr}_{\LocSys_\cG(X)}:\LocSys_\cG(X)\times X^I\to \LocSys_\cG(X).$$

\medskip

The assignment
$$(I,V_I,\CM_I)\mapsto \on{Ev}_{V_I,\CM_I}$$
defines a symmetric monoidal functor 
\begin{equation} \label{e:Loc}
\on{Ev}:\Rep(\cG,\Ran(X))\to \QCoh(\LocSys_\cG(X)).
\end{equation}

We have the following result:

\begin{prop}[D.G and J.~Lurie, unpublished]  \label{p:loc}
The functor $\on{Ev}$ admits a \emph{fully faithful} right adjoint.
\end{prop}

In other words, the above proposition says that $\QCoh(\LocSys_\cG(X))$ is a localization (a.k.a., Verdier quotient)
of $\Rep(\cG,\Ran(X))$ by a full subcategory, which is moreover a monoidal ideal.

\sssec{}  \label{sss:LocSys act}

According to \propref{p:loc}, given an action of the monoidal category $\Rep(\cG,\Ran(X))$ on some category $\bC$,
if this action factors through an action of $\QCoh(\LocSys_\cG(X))$ on $\bC$, then it does so uniquely. Moreover, this happens
if and only if the objects in
$$\on{ker}(\on{Ev})\subset \Rep(\cG,\Ran(X))$$
act on $\bC$ by zero.

\sssec{}

We are now ready to state the theorem (this is \cite[Theorem 4.5.2]{Ga3}) about the spectral decomposition of $\on{D}(\Bun_G(X))$ along $\LocSys_\cG(X)$.
Recall the action of $\Rep(\cG,\Ran(X))$ on $\on{D}(\Bun_G(X))$ from \secref{sss:Hecke act ult}.

\begin{thm}[V.~Drinfeld, D.G.] \label{t:gen vanishing}
The action of $$\on{ker}(\on{Ev})\subset \Rep(\cG,\Ran(X))$$ on $\on{D}(\Bun_G(X))$ is zero.
\end{thm}

According to \secref{sss:LocSys act}, from \thmref{t:gen vanishing} we obtain a canonically defined action of the monoidal category 
$\QCoh(\LocSys_\cG(X))$ on $\on{D}(\Bun_G(X))$, in such a way that the objects
$\on{Ev}_{V_I,\CM_I}$ (see \eqref{e:Ev}) acts as the endo-functors $\on{H}_{V_I,\CM_I}$. 

\medskip

We refer to this action as the `spectral decomposition of $\on{D}(\Bun_G(X))$ along $\LocSys_\cG(X)$'. 

\ssec{Relation to the `vanishing conjecture' of \cite{FGV}}

In the paper \cite{FGV} a certain conjecture was proposed (for which the sheaf-theoretic context can be either $\ell$-adic sheaves
or D-modules), and it was shown that this conjecture implies the existence of Hecke eigensheaves for $GL_n$. This conjecture
was subsequently proved in \cite{Ga1}.  

\medskip

In this subsection we will show that in the context of D-modules, the vanishing conjecture from
\cite{FGV} is a particular case of \thmref{t:gen vanishing} . 

\sssec{}

Let $G$ be $GL_n$. We consider the stack
$\Bun_n(X):=\Bun_{GL_n}$. For a non-negative integer $d$, let $\on{Mod}_{n,d}(X)$ be the stack classifying triples
$$(\CM,\CM',\alpha),$$
where $\CM,\CM'$ are rank-$n$ vector bundles on $X$, and $\alpha$ is an injection $\CM\hookrightarrow \CM'$
\emph{as coherent sheaves} so that the quotient $\CM'/\CM$ (which is a priori a torsion sheaf on $X$) has
length $d$. 

\medskip

We have the projections
$$\Bun_n(X)  \overset{\hl}\longleftarrow \on{Mod}_{n,d}(X) \overset{\hr}\longrightarrow \Bun_n(X)$$
where $\hl(\CM,\CM',\alpha)=\CM$ and $\hl(\CM,\CM',\alpha)=\CM'$. 

\medskip

Let 
$$\overset{\circ}{\on{Mod}}_{n,d}(X) \overset{j}\hookrightarrow \on{Mod}_{n,d}(X)$$ 
be the open substack corresponding to the condition that the quotient $\CM'/\CM$ be \emph{regular semi-simple}
(i.e., the direct sum of $d$ sky-scrapers concentrated in distinct points of $X$). 

\medskip

We have a projection
$$\overset{\circ}{\on{Mod}}_{n,d}(X)\overset{s}\to \overset{\circ}X{}^{(d)}$$
that remembers the support of $\CM'/\CM$, where $\overset{\circ}X{}^{(d)}\subset X^{(d)}$ is the open subscheme of
multiplicity-free divisors. 

\sssec{}

Let $E$ be a local system on $X$ (of an arbitrary finite rank). We form its symmetric power $E^{(d)}$, which is
a local system of rank $\on{rk}(E)\circ d$ when restricted to $\overset{\circ}X{}^{(d)}$. We consider the object
$$\fL_{E,d}\in \on{Mod}_{n,d}(X),$$
(known as Laumon's sheaf) defined as
$$\fL_{E,d}:=j_{!*}(s^*(E^{(d)})),$$
where $j_{!*}$ is the operation of Goresky-MacPherson extension (applied to the local system $s^*(E^{(d)})$ on
$\overset{\circ}{\on{Mod}}_{n,d}(X)$).

\medskip

We define the \emph{averaging functor}
$$\on{Av}_{E,d}:\on{D}(\Bun_n(X))\to \on{D}(\Bun_n(X)), \quad \on{Av}_{E,d}(\CF):=\hr_!(\hl{}^!(\CF)\sotimes \fL_{E,d}).$$

\medskip

The vanishing conjecture of \cite{FGV}/theorem of \cite{Ga1} says:

\begin{thm}  \label{t:vanishing}
Suppose that $E$ is irreducible and $\on{rk}(E)>n$ and $d>(2g-2)\cdot n\cdot \on{rk}(E)$. Then $\on{Av}_{E,d}=0$.
\end{thm}

\sssec{}

Let us specialize again to the case when $k$ is of characteristic $0$, and $\on{D}(-)$ is the derived category of D-modules.
We claim that in the case, \thmref{t:vanishing} is a tiny particular case of \thmref{t:gen vanishing}. 

\medskip

Indeed, it is easy to see that the functor $\on{Av}_{E,d}$ is given by the action of a particular object 
${\mathcal {Av}}_{E,d}\in \Rep(\cG,\Ran(X))$. Moreover, 
$$\on{Ev}({\mathcal {Av}}_{E,d})\in \QCoh(\LocSys_\cG(X))$$
is calculated as follows:

\medskip

We note that for $G=GL_n$, we have $\cG=GL_n$ and the fiber of $\on{Ev}({\mathcal {Av}}_{E,d})$ at $\sigma\in \LocSys_\cG(X)$
is given by
$$H(X^{(d)},(E\otimes E_\sigma)^{(d)}),$$
where $E_\sigma$ is the$n$-dimensional local system corresponding to $\sigma$. 

\medskip

Now, in order to deduce \thmref{t:vanishing}, we notice that the above cohomology identifies with
$$\on{Sym}^d(H(X,E\otimes E_\sigma)),$$
and the latter vanishes for all $\sigma$ under the conditions on $E$ and $d$ specified in the theorem.

\sssec{}

Let us also note that for $k=\BF_q$, \thmref{t:vanishing} says something quite non-trivial even about the classical Hecke operators
acting on $\on{Autom}(X)$.  Namely, it says that if $f$ is a joint eigenvector of the Hecke algebras $\BH(G)_x$ with characters
$$(\lambda_{x,1},...,\lambda_{x,n}/\text{permutation}),$$
then the Rankin-Selberg L-function 
$$L(E,f,t)=\underset{x}\Pi\, \frac{1}{1-t^{\on{deg}_x}\cdot (\lambda_{x,1}+...+\lambda_{x,n})\cdot \on{Tr}(\on{Frob}_x,E_x)}$$
is actually a polynomial of degree $\leq (2g-2)\cdot n\cdot \on{rk}(E)$. 

\section{Global unramified geometric Langlands} \label{s:global}

In this section we let the ground field $k$ be of characteristic zero. 

\ssec{Why the \emph{best hope} does not work} 

\sssec{}

The categorical global unramified geometric Langlands theory aimes to compare the categories
$\on{D}(\Bun_G(X))$ and $\QCoh(\LocSys_\cG(X))$. So far, by \thmref{t:gen vanishing}, we have that the monoidal category 
$\QCoh(\LocSys_\cG(X))$ acts on $\on{D}(\Bun_G(X))$. Therefore, the datum of a $\QCoh(\LocSys_\cG(X))$-linear functor
$$\QCoh(\LocSys_\cG(X))\to \on{D}(\Bun_G(X))$$ amounts to a choice of an object in $\on{D}(\Bun_G(X))$. 

\medskip

Based on many pieces of evidence, the object in $\on{D}(\Bun_G(X))$ that we want to choose for the \emph{global geometric Langlands}
equivalence is the `first Whittaker coefficient'\footnote{It is denoted $\on{Poinc}(\CW_{\on{vac}})$ in \cite[Sect. 5.7.4]{Ga3}.}. Thus,
we obtain a functor 
\begin{equation} \label{e:naive Langlands}
\BL_G:\QCoh(\LocSys_\cG(X))\to \on{D}(\Bun_G(X)).
\end{equation} 

\sssec{}  \label{sss:arth}

When $G=T$ is a torus, the above functor $\BL_G$ is the generalized Fourier-Mukai transform studied by G.~Laumon. In particular,
it is an equivalence.

\medskip

However, the functor $\BL_G$ cannot be an equivalence as long $G$ is non-commutative. The reason is that there are objects in 
$\on{D}(\Bun_G(X))$ that are Whittaker-degenerate. For example, the constant D-module on $\Bun_G(X)$ is Whittaker-degenerate. 

\begin{rem}
Another heuristic piece of evidence for why the category $\on{D}(\Bun_G(X))$ cannot be equivalent to $\QCoh(\LocSys_\cG(X))$ comes
from the classical theory of automorphic functions: 

\medskip

It is known that automorphic representations are parameterized not by Langlands
parameters (i.e., Galois representations) but by \emph{Arthur parameters}, where the latter are conjugacy lasses of pairs $(\sigma,A)$ 
with $\sigma$ being a representation of the Galois group of $X$ into $\cG$, and $A$ is a nilpotent element of the Lie algebra of $\cG$ 
that commutes with $\sigma$.
\end{rem} 

\ssec{How to make the \emph{best hope} work?}
 
In \cite{AG1} an idea was suggested as to how one can modify the \emph{best hope} to make it work. 

\medskip

This modification
consists of tweaking the Galois side, i.e., replacing $\QCoh(\LocSys_\cG(X))$ by some other (but closely related)
category, while leaving the automorphic side intact. This tweak happens within homological algebra and has to do
with the fact that the category $\QCoh(\LocSys_\cG(X))$ is of infinite cohomological dimension. 

\sssec{}

In order to explain it we consider the following example. Consider the differential graded algebra $A:=k[\epsilon]$, where $\epsilon$ is
a free generator in degree $-1$ and its differential is zero (so the differential on all of $A$ is actually zero).

\medskip

Consider the (derived) category $A\mod$ of $A$-modules.  Inside we consider the full subcategory
$$A\mod^{\on{perf}}\subset A\mod$$
spanned by \emph{perfect complexes}, i.e., those objects that can be obtained by a finite process of
taking directs sums, summands and cones from the object $A\in A\mod$ itself. 

\medskip

We can also consider the full subcategory $A\mod^{\on{f.g.}}\subset A\mod$ spanned by objects that have
finite-dimensional cohomologies, all of which are finite-dimensional as vector spaces over $k$. 

\medskip

Since $A\in A\mod^{\on{f.g.}}$, we have the inclusion 
$$A\mod^{\on{perf}}\subset A\mod^{\on{f.g.}},$$
but it is not an equality. In fact, the Verdier quotient $A\mod^{\on{f.g.}}/A\mod^{\on{perf}}$ is equivalent to the
category $B\mod^{\on{f.g.}}$, where $B=k[\eta,\eta^{-1}]$, where $\eta$ is a free generator in degree $2$
and its differential is zero. 

\begin{rem}
A similar phenomenon in the category of representations of a finite group with torsion coefficients leads to
the notion of \emph{Tate cohomology}.
\end{rem}

\sssec{}

More generally, let $V$ be a finite-dimensional vector space, and consider the DG algebra (with zero differential)
$A:=\Sym(V[1])$. As above we consider the categories
$$A\mod^{\on{perf}}\subset A\mod^{\on{f.g.}}.$$

However, one can now notice that to every conical Zariski-closed subset $\CN\subset V$ one can attach a full subcategory
$$A\mod^{\on{f.g.}}_\CN\subset A\mod^{\on{f.g.}},$$
such that:
$$
\begin{cases}
&A\mod^{\on{f.g.}}_\CN=A\mod^{\on{perf}} \text{ if } \CN=0; \\
&A\mod^{\on{f.g.}}_\CN=A\mod^{\on{f.g.}} \text{ if } \CN=V.
\end{cases}
$$
 
\sssec{}

Even more generally, let $Y$ be a scheme (or algebraic stack), which is a \emph{locally complete intersection}. To
$Y$ one attaches another scheme/stack (see \cite[Sect. 2.3]{AG1}), denoted $\on{Sing}(Y)$, whose $k$-points are pairs $(y,\xi)$, 
where $y$ is a $k$-point of $Y$, and $\xi$ is an element of the vector space $H^{-1}(T^*_y(Y))$, where $T^*_y(Y)$
is the \emph{derived cotangent space} at $y$, i.e., the fiber of the cotangent complex\footnote{We recall that 
locally complete intersections are characterized by the property that $H^{-i}(T^*_y(Y))=0$ for all $y$ and $i>1$.}
of $Y$ at $y$. 
 
\medskip 

Let $\Coh(Y)$ be the full subcategory of $\QCoh(Y)$ consisting of objects with finitely many cohomologies, each of which
is coherent (i.e., locally finitely generated) as a sheaf on $Y$. In \cite[Sect. 4]{AG1} the following construction is performed: to every conical 
Zariski-closed subset $\CN\subset \on{Sing}(Y)$ one attaches a full subcategory
$$\Coh_\CN(Y)\subset \Coh(Y).$$

Again, we have:
$$
\begin{cases}
&\Coh_\CN(Y)=\on{Perf}(Y) \text{ if } \CN=\{\underset{y\in Y}\cup\, {(y,0)}\}; \\
&\Coh_\CN(Y)=\Coh(Y)  \text{ if } \CN=\Sing(Y),
\end{cases}
$$
where $\on{Perf}(Y)\subset \Coh(Y)$ is the subcategory of perfect objects (complexes that locally on $Y$ can be represented by a 
finite complex of free sheaves of finite rank). 

\sssec{}

The enlargement 
$$\on{Perf}(Y)\rightsquigarrow \Coh_\CN(Y)$$ is
exactly the tweak that we will perform on the Galois side of the global unramified geometric Langlands theory. However, there is one point of difference.

\medskip

For multiple reasons, it is more convenient to work with large (technical term: \emph{cocomplete compactly generated}) categories 
(such as $\QCoh(Y)$), i.e., categories that admit arbitrary direct sums (and generated by a set of compact objects).
The datum of such a category is equivalent to the datum of
its full subcategory of compact objects (in the case of $\QCoh(Y)$, its subcategory of compact objects is exactly $\on{Perf}(Y)$),
which is a small category. The inverse procedure (recovering a large category from a small one) is called \emph{ind-completion},
see \cite[Sect. 7.2]{GR1}.

\medskip

The large subcategory corresponding to $\Coh_\CN(Y)$ is denoted $\IndCoh_\CN(Y)$. The inclusion 
$\on{Perf}(Y)\subset \Coh_\CN(Y)$ extends to a fully faithful functor
$$\QCoh(Y)\hookrightarrow \IndCoh_\CN(Y),$$
which admits a right adjoint, given by ind-extending the tautological embedding
$$\Coh_\CN(Y)\hookrightarrow \Coh(Y)\hookrightarrow \QCoh(Y).$$

\sssec{}  \label{sss:t-str}

For any $\CN$, the category $\IndCoh_\CN(Y)$ carries a t-structure\footnote{This is one of the reasons to work 
with the large category $\IndCoh_\CN(Y)$ as opposed to the small category $\Coh_\CN(Y)$.}, and the functor
$$\IndCoh_\CN(Y)\to \QCoh(Y)$$
(right adjoint to the tautological inclusion) 
is t-exact. Moreover, the above functor \emph{induces an equivalence of the corresponding bounded below
subcategories}
$$\IndCoh_\CN(Y)^+\to \QCoh(Y)^+,$$
see \cite[Sect. 4.4]{AG1}. 

\medskip

So, the difference between $\QCoh(Y)$ and $\IndCoh_\CN(Y)$ `occurs at $-\infty$'. Note that there is no
contradiction here: the t-structure on $\IndCoh_\CN(Y)$ is non-separated, that is, there are non-zero objects all of whose
cohomologies vanish. 

\ssec{Back to $\LocSys_\cG(X)$}

The modification of the Galois side of geometric Langlands , proposed in \cite{AG1}, is of the form
$$\IndCoh_\CN(\LocSys_\cG(X)),$$
for a particular conical Zariski-closed subset $\CN\subset \Sing(\LocSys_\cG(X))$.

\sssec{}

First, we describe the stack $\Sing(\LocSys_\cG(X))$. By unwinding the definition of the cotangent complex (see \cite[Sect. 10.4.6]{AG1}),
we obtain that $\Sing(\LocSys_\cG(X))$ is the moduli stack of pairs $(\sigma,A)$, where $\sigma\in \LocSys_\cG(X)$ and $A$ is a \emph{horizontal}
section of the local system associated with the co-adjoint representation of $\cG$.

\medskip

Choosing an ad-invariant symmetric bilinear form on $\cg$, we can think of $A$ as a section of the local system associated with the 
adjoint representation of $\cG$. We let
$$\CN\subset \Sing(\LocSys_\cG(X))$$
be the \emph{global nilpotent cone}, i.e., the set of those $(\sigma,A)$ for which $A$ is nilpotent as a section of the local system of
Lie algebras $\cg_\sigma$ (equivalently, the value of $A$ in the fiber of $\cg_\sigma$ at some/every point of $X$ 
should be nilpotent).  

\sssec{}

Thus, the proposed category on the Galois side of global unramified geometric Langlands is
$\IndCoh_\CN(\LocSys_\cG(X))$ for the above choice of $\CN$.

\medskip

We note that if $G=T$ is a torus, the nilpotent cone in $\cg$ is zero. So, in this case
$$\IndCoh_\CN(\LocSys_\cT(X))=\QCoh(\LocSys_\cT(X)),$$
i.e., the Galois side is the same as in the original \emph{best hope} (as it should be, because the \emph{best hope}
is realized by the Fourier-Mukai transform). 

\medskip

However, this modification is nontrivial as soon as $G$ is non-commutative. The \emph{most singular} point of
$\LocSys_\cG(X)$ is one corresponding to the \emph{trivial local system}. Around this point, the difference between
$\IndCoh_\CN(\LocSys_\cG(X))$ and the initial $\QCoh(\LocSys_\cG(X))$ is the largest.

\medskip

Consider now the open substack
$$\LocSys^{\on{irred}}_\cG(X)\subset \LocSys_\cG(X)$$
consisting of irreducible local systems. It is easy to see that the corresponding inclusion
$$\QCoh(\LocSys^{\on{irred}}_\cG(X))\subset \IndCoh_\CN(\LocSys^{\on{irred}}_\cG(X))$$
is an equality. So, the modification does not affect the irreducible locus\footnote{This fact may arouse suspicions
in the validity of the proposed form of the geometric Langlands conjecture for groups other than $GL_n$.}.
 
\sssec{}

Thus, the proposed version of the global unramified geometric Langlands equivalence reads as follows:

\begin{conj} \label{c:geom Langlands}
There is a canonically defined equivalence of categories
$$\BL_G:\IndCoh_\CN(\LocSys_\cG(X))\to \on{D}(\Bun_G(X)).$$
\end{conj}

The above statement of \conjref{c:geom Langlands} is too loose. The paper \cite{Ga3} lists a list of 
compatibility requirements that fix $\BL_G$ is uniquely. 

\begin{rem}
One can view \conjref{c:geom Langlands} as restoring the Arthur parameters (as opposed to just Langlands
parameters) that were missing in the original  \emph{best hope}. Indeed, they appear as obstructions to 
\emph{temperedness}, i.e., as obstructions for an object of $\on{D}(\Bun_G(X))$ to be in the essential
image of $\QCoh(\LocSys_\cG(X))$.
\end{rem} 

\begin{rem}
Recall that in \secref{sss:t-str} we said that the difference between the categories $\IndCoh_\CN(\LocSys_\cG(X))$
and $\QCoh(\LocSys_\cG(X))$ `occurs at $-\infty$' with respect to their respective t-structures.  On the other hand, the
failure of the functor \eqref{e:naive Langlands} to be an equivalence happens already at the level of the
corresponding bounded categories: indeed, recall that the constant D-module on $\Bun_G(X)$ is not in the image of
$\QCoh(\LocSys_\cG(X))$. This `contradiction' is explained by the fact that the functor $\BL_G$ in \conjref{c:geom Langlands} is 
of inifinite cohomological amplitude. 
\end{rem}

\ssec{What is known?}

An outline of how one might go about proving \conjref{c:geom Langlands} was proposed in \cite{Ga3}. Here will summarize the main ideas
of the proposed proof and comment on its status. It consists of the following steps.

\sssec{}

On the automorphic side one constructs a category, denoted $\Whit^{\on{ext}}_{G,G}(X)$ (called the \emph{extended Whittaker category},
see \cite[Sect. 8.2]{Ga3}) and a functor of Whittaker expansion
\begin{equation} \label{e:coeff functor}
\on{coeff}^{\on{ext}}_{G,G}:\on{D}(\Bun_G(X))\to \Whit^{\on{ext}}_{G,G}(X).
\end{equation} 

\medskip

In \cite[Conjecture 8.2.9]{Ga3} it is conjectured that the functor $\on{coeff}^{\on{ext}}_{G,G}$ is fully faithful, and this conjecture
is proved in \cite{Ber} for $G=GL_n$. 

\sssec{}

The category $\Whit^{\on{ext}}_{G,G}(X)$ can be thought of as fibered over the space of characters $\fch$ 
of $N(\BA)$ trivial on $N(K)$, where for each $\chi\in \fch$ we consider the category of D-modules on $G(\BO)\backslash G(\BA)$
that transform according to $\chi$ with respect to the action of $N(\BA)$. 

\medskip

The space $\fch$ splits according to the pattern of how degenerate the character is, i.e., it is a union of locally closed
subspaces $\fch^P$, where $P$ runs over the poset of standard parabolics of $G$.  For each $P$, we obtain the corresponding
full subcategory $\Whit_{G,P}(X)$.

\medskip

For example, for $P=G$, we obtain the the usual Whittaker category, denoted $\Whit_{G,G}(X)$, see \cite[Sect. 5]{Ga3}.

\medskip

For $P=B$, the category $\on{Whit}_{G,B}(X)$ is the \emph{principal series category} of \cite[Sect. 6]{Ga3}. 

\sssec{}

On the Galois side one constructs a category, denoted $\on{Glue}(\cG)_{\on{spec}}$, and a functor 
\begin{equation} \label{e:spec functor}
\IndCoh_\CN(\LocSys_\cG(X))\to \on{Glue}(\cG)_{\on{spec}}.
\end{equation} 

In \cite{AG2} it is shown that \eqref{e:spec functor} is fully faithful. 

\sssec{}

The category $\on{Glue}(\cG)_{\on{spec}}$ is glued from the categories 
$$\QCoh_{\on{conn}/\LocSys_\cG(X)}(\LocSys_\cP(X)),$$
where $\cP$ runs through the poset of standard parabolics of $\cG$. 

\medskip

In the above formula, $\QCoh_{\on{conn}/\LocSys_\cG(X)}(\LocSys_\cP(X))$ is the (derived) category of quasi-coherent sheaves on the stack $\LocSys_\cP(X)$, 
endowed with a \emph{connection} along the fibers of the map $\LocSys_\cP(X)\to \LocSys_\cG(X)$. These notions need to be understood in the
sense of derived algebraic geometry, see \cite[Sect. 6.5]{Ga3}. 

\medskip

For example, for $P=G$, we have $$\QCoh_{\on{conn}/\LocSys_\cG(X)}(\LocSys_\cP(X))=\QCoh(\LocSys_\cG(X)),$$ i.e., this is the usual
(unmodified) derived category of quasi-coherent sheaves on $\LocSys_\cG(X)$.

\sssec{}

Assuming \conjref{c:geom Langlands} for \emph{proper Levi subgroups} of $G$, and certain auxiliary results,  
one constructs a \emph{fully faithful} functor 
\begin{equation} \label{e:loc functor}
\on{Glue}(\cG)_{\on{spec}}\to \Whit^{\on{ext}}_{G,G}(X).
\end{equation} 

\medskip

The construction of \eqref{e:loc functor} with the required properties is complete for $G=GL_2$, and it is a question of time
before it becomes available for any $G$.

\sssec{}

The functor \eqref{e:loc functor} is glued from the fully faithful functors
$$\QCoh_{\on{conn}/\LocSys_\cG(X)}(\LocSys_\cP(X))\to \Whit_{G,P}(X).$$

For example, for $P=G$, the corresponding functor 
$$\QCoh(\LocSys_\cG(X))\to \Whit_G(X)$$
is the composition of the functor
$$\QCoh(\LocSys_\cG(X))\to \Rep(\cG,\Ran(X)),$$
right adjoint to the functor $\on{Ev}$ of \eqref{e:Loc} followed by the Casselman-Shalika equivalence
$$\Rep(\cG,\Ran(X))\simeq \Whit_G(X).$$

Here, $\Whit_G(X)$ is a slight modification of $\Whit_{G,G}(X)$ that has to do with the center of $G$,
see \cite[Sect. 5.6.7]{Ga3}. 

\sssec{}

Let us now assume that the functor \eqref{e:coeff functor} is fully faithful (which we know for $GL_n$), and that the functor
\eqref{e:loc functor} with the required properties exists. Let us see how this helps to prove \conjref{c:geom Langlands}. 

\medskip

Consider the composition of the functors \eqref{e:spec functor} and \eqref{e:loc functor}, which is a fully faithful functor
\begin{equation} \label{e:comp} 
\IndCoh_\CN(\LocSys_\cG(X))\to \Whit^{\on{ext}}_{G,G}(X).
\end{equation}

In \cite[Sect. 10]{Ga3}, one exhibits a collection of objects $\CF_\alpha\in \on{D}(\Bun_G(X))$ and a collection of objects
$\CM_\alpha\in \IndCoh_\CN(\LocSys_\cG(X))$ such that for every $\alpha$, the image of $\CF_\alpha$ under \eqref{e:coeff functor}
is isomorphic to the image of $\CM_\alpha$ under \eqref{e:comp}.

\medskip

In addition, it is shown that the objects $\CF_\alpha$ \emph{generate} $\on{D}(\Bun_G(X))$. And is conjectured (and established for $G=GL_n$)
that the objects $\CM_\alpha$ generate $\IndCoh_\CN(\LocSys_\cG(X))$. This implies that the essential images of \eqref{e:coeff functor} and
\eqref{e:comp} in $\Whit^{\on{ext}}_{G,G}(X)$, being generated by the same collection of objects, coincide, thus providing the sought-for equivalence $\BL_G$. 

\begin{rem}
While all the preceding steps in the proposed proof of \conjref{c:geom Langlands} were geometric in nature (i.e., used the standard sheaf-theoretic
functors on the categories D-modules when working on the automorphic side) and had clear counterparts in the classical theory of automorphic 
functions, the construction of the objects $\CF_\alpha$ (resp., $\CM_\alpha$) is different in nature, and is based on the ideas from \cite{BD}:

\medskip

On the automorphic side, the objects $\CF_\alpha$ are obtained by the \emph{localization functor} from modules over the Kac-Moody algebra
at the critical level. On the Galois side, the objects $\CM_\alpha$ are obtained by taking direct images along a map to $\LocSys_\cG(X)$
from the scheme classifying $\cG$-\emph{opers} on $X$ (see \cite[Sect. 10]{Ga3}).

\end{rem}

\section{Local geometric Langlands}  \label{s:loc}

\ssec{What is the object of study on the representation-theoretic side?}

\sssec{}

Recall that in the classical global Langlands theory, the object of study on the  representation-theoretic (=automorphic) side
is the space of functions on the quotient $G(\BA)/G(K)$, viewed
as a representation of the group $G(\BA)$. In the unramified case, the corresponding object of study is the space of functions on 
$G(\BO)\backslash G(\BA)/G(K)$, viewed as a module over the Hecke algebra.

\medskip

By contrast, the object of study on the representation-theoretic side in the classical \emph{local} theory is the \emph{category} of
representations of the group $G(\CK)$, where $\CK$ is a local field. So, by going from global to local we raise the level by one in the
hierarchy 
$$\on{Elements\, of \,a\, Set}\to \on{Objects\, of\, a\, Category}\to \on{Objects \,in\, a\, 2-Category}.$$

\medskip

In the global geometric Langlands theory, in the unramified case, the object of study on the representation-theoretic (=automorphic) side
was the \emph{category} $\on{D}(\Bun_G(X))$, viewed
as acted on by the Hecke functors. 

\medskip

Hence, by the above analogy, on the representation-theoretic in the local geometric theory, 
the object of study should be a certain 2-category, attached to the group $G$ and the local field $\CK=k\ppart$.  

\medskip

We stipulate that the 2-category in question is that of \emph{categories equipped with an action of $G(\CK)$}.  We will now explain what we mean
by this.

\sssec{}

First, when we say `category' in the above context, we mean a $k$-linear \emph{DG category}\footnote{All our DG categories are assumed \emph{cocomplete}.},
defined, e.g., as in \cite[Sect. 10]{GR1}. The important fact is that the totality of such categories and $k$-linear functors between 
them\footnote{All our functors are assumed \emph{continuous}, i.e., preserving infinite direct sums} forms an $(\infty,2)$-category, denoted $\on{DGCat}$, 
equipped with a symmetric monoidal structure, called the \emph{Lurie tensor product}. 

\medskip

Second, when we write $G(\CK)$ we mean the group ind-scheme, defined as a functor on the category of affine schemes by
$$\Hom(\Spec(A),G(\CK)):=\Hom(\Spec(A\ppart,G).$$

\medskip

We now have to define the notion of action\footnote{On the geometric/automorphic side of Langlands, when we talk about 
actions of groups on categories, we mean \emph{strong} actions.} of $G(\CK)$ on a category. The corresponding general 
notion has been developed in \cite{Ga4}. 

\medskip

However, one can give also the following explicit definition: according to \cite{Ber}, we have a well-defined category 
$\on{D}(G(\CK))$ of D-modules on $G(\CK)$; the group structure on $G$ defines on $\on{D}(G(\CK))$ a monoidal 
structure. I.e., $\on{D}(G(\CK))$ acquires a structure of \emph{associative algebra} in the monoidal category
$\on{DGCat}$. 

\medskip

We define the notion of category equipped with an action of $G(\CK)$ to be a module over
$\on{D}(G(\CK))$ in $\on{DGCat}$. The totality of such has a structure of $(\infty,2)$-category (see \cite[Sect. 8.3]{GR1});
we denote it by $G(\CK)\mmod$. 

\sssec{}

Here are some examples of objects of $G(\CK)\mmod$. 

\medskip

\noindent(i) The first example is $\bC:=\on{D}(G(\CK))$, equipped with an action on itself by left multiplication.

\medskip

\noindent(ii) For any subgroup $H\subset G(\CK)$, we can take $\bC:=\on{D}(G(\CK)/H)$. 

\medskip

As particular cases of the above example (ii), we can take $H=G(\CO)$ or $H=I$, the latter
being the Iwahori subgroup. The resulting categories are the categories of D-modules on the
affine Grassmannian $\on{Gr}_G$ and the affine flag scheme $\on{Fl}_G$, respectively. 

\medskip

\noindent(iii) We can take $\bC=\hg_\kappa\mod$, i.e., the category of representations of
the Kac-Moody algebra for any \emph{integral} level $\kappa$ (see \cite[Sect. 23]{FG} for
the definition); here the action of $G(\CK)$ comes from its adjoint action on $\hg_\kappa$. 

\medskip

\noindent(iv) We consider the stack $\Bun_G(X)^{\on{level}_x}$, classifying principal $G$-bundles
on the curve $X$ equipped with a full level structure at a point $x$ (i.e., a trivialization over
the formal neighborhood of $x$). We take $\bC:=\on{D}(\Bun_G(X)^{\on{level}_x})$. This
is the object of $G(\CK)\mmod$, corresponding to the global geometric Langlands theory
with ramification allowed at $x$. 

\ssec{The object of study on the Galois side}
 
\sssec{}

Recall that in the global unramified geometric theory, the object of study on the Galois side was (a modification of) the derived category of 
quasi-coherent sheaves on the stack $\LocSys_\cG(X)$ that classifies $\cG$-local systems on the curve $X$.

\medskip

Based on the analogy with the local classical theory, the object of study on the Galois side in the \emph{local geometric theory}
should be a certain 2-category attached to the space $\LocSys_\cG(\oD)$ of $\cG$-local systems on the \emph{puntured formal disc}
$\oD$. 

\medskip

In what follows we will explain what we mean by the space $\LocSys_\cG(\oD)$, and what the resulting 2-category
(in fact, $(\infty,2)$-category) is.

\sssec{}  \label{sss:loc Galois side}

Recall also that in the global case, the category that one obtains from $\LocSys_\cG(X)$ in the most tautological way, i.e.,
$\QCoh(\LocSys_\cG(X))$, did not quite match the automorphic side--we needed to introduce a correction that had to
do with the difference between perfect complexes and coherent ones. 

\medskip

The $(\infty,2)$-category $\on{ShvCat}(\LocSys_\cG(\oD))$ that we will define below is the counterpart of $\QCoh(\LocSys_\cG(X))$.
It will be responsible for the tempered part of $G(\CK)\mmod$. 

\medskip

The extension of $\on{ShvCat}(\LocSys_\cG(\oD))$ that takes into account all local Arthur parameters (as opposed to just Langlands
parameters) has been recently proposed by D.~Arinkin. But we will not explicitly discuss it in this talk. 

\ssec{Sheaves of categories}

In order to talk about sheaves of categories, we need to place ourselves in the context of \emph{derived algebraic geometry}.
Thus, in what follows, when we say `affine scheme' we shall mean a derived affine scheme over $k$. By definition, the category
of such is the opposite of the category of \emph{connective}\footnote{Connective=concentrated in non-positive cohomological
degrees.} commutative DG algebras over $k$, see \cite[Sect. 1.1]{GR2}.

\sssec{}

For an affine scheme $S$, we consider the (symmetric) monoidal category $\QCoh(S)$. This is a (commutative) algebra
object in the symmetric monoidal category $\on{DGCat}$. 

\medskip

We let $\on{ShvCat}(S)$ denote the $(\infty,2)$-category $\QCoh(S)\mmod$, i.e., that of $\QCoh(S)$-modules in
the (symmetric) monoidal $(\infty,2)$-category $\on{DGCat}$. 

\medskip

The assignment $S\mapsto \on{ShvCat}(S)$ is a functor from 
$(\affdgSch)^{\on{op}}$ to the $\infty$-category of $(\infty,2)$-categories.

\sssec{}

Let now $\CY$ be an arbitrary \emph{prestack}, i.e., a functor
$$(\affdgSch)^{\on{op}}\to \Spc,$$
where $\Spc$ is the $\infty$-category of \emph{spaces}\footnote{Space=$\infty$-groupoid.}. 

\medskip

We define the $(\infty,2)$-category $\on{ShvCat}(\CY)$ to be
$$\underset{S\overset{y}\to \CY}{\on{lim}}\, \on{ShvCat}(S),$$
where the index category is $((\affdgSch)_{/\CY})^{\on{op}}$ and the limit
is taken in the $\infty$-category of $(\infty,2)$-categories.

\sssec{}

In other words, informally, an object of $\on{ShvCat}(\CY)$ is an assignment for every 
\begin{equation} \label{e:eval sheaf}
(S\overset{y}\to \CY)\, \rightsquigarrow\, \bC_{S,y}\in \QCoh(S)\mmod,
\end{equation} 
and for every 
$$((S_1,y_1),(S_2,y_2),S_1\overset{f}\to S_2, y_2\circ f\sim y_1)\,  \rightsquigarrow\,
\QCoh(S_1)\underset{\QCoh(S_2)}\otimes \bC_{S_2,y_2} \simeq \bC_{S_1,y_1}.$$
This assignment must be endowed with a \emph{homotopy-coherent system of compatibilities} for compositions. 

\medskip

We call objects of $\on{ShvCat}(\CY)$ `sheaves of categories over $\CY$'.  

\sssec{}

The most basic
example of a sheaf of categories is $\QCoh_{/\CY}$ (not to be confused with the \emph{category}
$\QCoh(\CY)$, discussed below). Namely, in terms of the assignment \eqref{e:eval sheaf}, 
the object $\QCoh_{/\CY}$ assigns to 
$$(S,y)\rightsquigarrow   \QCoh(S)\in \QCoh(S)\mmod.$$

\sssec{}

For $\CY$ as above we can consider the (symmetric) monoidal category
\begin{equation} \label{e:qcoh}
\QCoh(\CY):=\underset{S\overset{y}\to \CY}{\on{lim}}\, \QCoh(S).
\end{equation} 

This is what one calls the \emph{(derived) category of quasi-coherent sheaves on a prestack}. (In the case of $\CY=\LocSys_\cG(X)$,
this is the category $\QCoh(\LocSys_\cG(X))$ considered in the previous sections.) 

\medskip

In other words, informally, an object of $\QCoh(\CY)$ is an assignment for every 
\begin{equation} \label{e:eval sheaf q}
(S\overset{y}\to \CY)\, \rightsquigarrow\, \CF_{S,y}\in \QCoh(S),
\end{equation} 
and for every 
$$((S_1,y_1),(S_2,y_2),S_1\overset{f}\to S_2, y_2\circ f\sim y_1)\,  \rightsquigarrow\,
f^*(\CF_{S_2,y_2}) \simeq \CF_{S_1,y_1}.$$
This assignment must be endowed with a \emph{homotopy-coherent system of compatibilities} for compositions. 

\sssec{}

Here is another candidate for what we might call a `sheaf of categories': we can consider the $(\infty,2)$-category 
$$\QCoh(\CY)\mmod.$$

\medskip

Note that if $\CY$ is (representable by) an affine scheme $S$, we have a tautological equivalence
$$\QCoh(S)\mmod\simeq \on{ShvCat}(S).$$

\sssec{}

For a general prestack $\CY$, the above two $(\infty,2)$-categories are related by a pair of adjoint functors
\begin{equation} \label{e:1-aff}
\bLoc:\QCoh(\CY)\mmod\rightleftarrows \on{ShvCat}(\CY):\bGamma.
\end{equation}

In terms of the assignment \eqref{e:eval sheaf}, the functor $\bGamma$ sends a sheaf of categories to
$$\underset{S\overset{y}\to \CY}{\on{lim}}\, \bC_{S,y}\in \on{DGCat},$$
which is equipped with a natural action of \eqref{e:qcoh}.

\medskip

The functor $\bLoc$ sends $\bC\in \QCoh(\CY)\mmod$ to
$$(S,y)\, \rightsquigarrow\, \QCoh(S)\underset{\QCoh(\CY)}\otimes \bC.$$

\medskip

We shall say that a prestack $\CY$ is \emph{1-affine} if the functors \eqref{e:1-aff} are mutually inverse equivalences.

\medskip

Tautologically, every affine scheme is 1-affine. 

\sssec{}  \label{sss:ex 1-aff}

Here are some examples of prestacks that are (or are not) 1-affine (these examples are taken from 
\cite[Sect. 2]{Ga4}):

\medskip

\noindent(i) Any quasi-compact, quasi-separated algebraic space (in particular, scheme) is 1-affine.

\medskip

\noindent(ii) Any quasi-compact algebraic stack of finite type with an affine diagonal is 1-affine.

\medskip

\noindent(iii) For a non-trivial connected algebraic group $G$, the quotient $\on{pt}/G(\CO)$ is \emph{not} 1-affine.
(This is not in contradiction with example (ii), because the finite-type condition is violated.) 

\medskip

\noindent(iv) The ind-scheme $\BA^\infty=\underset{i}{\on{colim}}\, \BA^i$ is \emph{not} 1-affine.

\medskip

In general, one can say that infinite-dimensionality is typically an obstruction to 1-affineness.

\ssec{The space of local systems on the formal punctured disc}

\sssec{}

We now introduce the space $\LocSys_\cG(\oD)$, which is the main geometric player on the Galois side of
the local geometric Langlands.

\medskip

We start with the space of $\cg$-valued connection forms on $\oD$, i.e., $\cg\otimes \omega_\CK$. This is
an ind-scheme (of infinite type). A choice of a uniformizer in $\CK$ identifies $\cg\otimes \omega_\CK$ with $\cg(\CK)$. 

\medskip

Now, the group $G(\CK)$ acts on $\cg\otimes \omega_\CK$ by \emph{gauge transformations}.

\medskip

We define $\LocSys_\cG(\oD)$ to be the \emph{prestack quotient} $\cg\otimes \omega_\CK/G(\CK)$. 

\sssec{}  \label{sss:loc expect}

As was mentioned in \secref{sss:loc Galois side}, the $(\infty,2)$-category
$$\on{ShvCat}(\LocSys_\cG(\oD))$$
plays the same role vis-\`a-vis $G(\CK)\mmod$ as $\QCoh(\LocSys_\cG(X))$ did vis-\`a-vis $\on{D}(\Bun_G)$.

\medskip

In particular, we expect that 

\medskip

\noindent(i) The $(\infty,2)$-category $\on{ShvCat}(\LocSys_\cG(\oD))$ (which is equipped with a natural (symmetric) monoidal structure)
acts on $G(\CK)\mmod$;

\medskip

\noindent(ii) $\on{ShvCat}(\LocSys_\cG(\oD))$ is equivalent to the full subcategory of $G(\CK)\mmod$,
consisting of tempered objects. 

\sssec{}

We propose:

\begin{conj}  \label{c:1-aff} 
The prestack $\LocSys_\cG(\oD)$ is 1-affine.
\end{conj}

This conjecture would imply that the two possible candidates for the notion of category over $\LocSys_\cG(\oD)$, namely,
$$\on{ShvCat}(\LocSys_\cG(\oD)) \text{ and } \QCoh(\LocSys_\cG(\oD))\mmod,$$
are equivalent.

\medskip

\conjref{c:1-aff} is a theorem when $G$ is a torus.

\begin{rem} 
Let us point out that the validity of \thmref{t:1-aff}(a) is something one should not take for granted. 
Indeed,  $\LocSys_\cG(\oD)$ is obtained
by quotienting $\cg\otimes \omega_\CK$ by $G(\CK)$. Now, the ind-scheme $\cg\otimes \omega_\CK$ \emph{is not} 1-affine
(because it contains Example(iv) from \secref{sss:ex 1-aff}). Furthermore, quotienting by groups of infinite type also
usually destroys 1-affineness (see Example(iii) from \secref{sss:ex 1-aff}).  Thus, for example, if instead of of the \emph{gauge} action of 
$G(\CK)$ on $\cg\otimes \omega_\CK$ we considered the adjoint action, the resulting quotient would \emph{not} be 1-affine. 

\medskip 

However, we expect that $\LocSys_\cG(\oD)$ still manages to be 1-affine: the infinite-dimensional ind-direction in $\cg\otimes \omega_\CK$ 
that prevents it from being 1-affine gets `eaten up' by the ind-direction in $G(\CK)$, and the action of $G(\CO)$ is 
`free modulo something finite-dimensional', i.e., the example Example(iii) from \secref{sss:ex 1-aff} does not occur. 
So, $\LocSys_\cG(\oD)$ does not have infinite-dimensional features that prevent it from being 1-affine, 
and yet it is \emph{not} locally of finite type (except if $G$ is a torus).

\end{rem}

\sssec{}

The following partial result towards \conjref{c:1-aff} has been recently established in \cite{Ras}:

\begin{thm}[S.~Raskin] \label{t:1-aff} \hfill

\smallskip

\noindent{\em(a)} The functor $\bLoc$ for $\LocSys_\cG(\oD)$ is fully faithful.

\smallskip

\noindent{\em(b)} The category $\QCoh(\LocSys_\cG(\oD))$ is compactly generated.

\end{thm}

\end{document}